\documentclass[11pt]{amsart}
\usepackage[latin1]{inputenc}
\usepackage{amsmath,amsthm, amscd, amssymb, amsfonts}
\usepackage[dvips]{graphicx}
\numberwithin{equation}{section}
\newcommand{\s}{\mathbb{S}}

\theoremstyle{plain}

\newtheorem{teor}{Theorem}[section]
\newtheorem{corol}[teor]{Corollary}
\newtheorem{prop}[teor]{Proposition}
\newtheorem{lem}[teor]{Lemma}
\newtheorem{afir}[teor]{Claim}

\theoremstyle{definition}
\newtheorem{defin}[teor]{Definition}

\newtheorem{ejem}[teor]{Example}

\theoremstyle{remark}
\newtheorem{obs}[teor]{Remark}

\newcommand\ad{\operatorname{ad}}

\newcommand\cop{\operatorname{cop}}

\newcommand\Rep{\operatorname{Rep}}
\newcommand\Ind{\operatorname{Ind}}

\newcommand\vect{\operatorname{Vec}}

\begin{document}

\title[Simple deformations of finite groups]{Simple Hopf algebras and deformations of finite groups}
\author{C\' esar Galindo and Sonia Natale}
\address{Facultad de Matem\'atica, Astronom\'\i a y F\'\i sica,
Universidad Nacional de C\'ordoba, CIEM -- CONICET, (5000) Ciudad
Universitaria, C\'ordoba, Argentina}
\email{galindo@mate.uncor.edu, natale@mate.uncor.edu}
\thanks{This work was partially supported by CONICET,
CONICOR, Fundaci\' on Antorchas  and Secyt (UNC)} 
\date{September 5, 2006}

\begin{abstract}
We show that certain twisting deformations of a family of
supersolvable groups are simple as Hopf algebras. These groups are
direct products of two generalized dihedral groups. Examples of
this construction arise in dimensions $60$ and $p^2q^2$, for prime
numbers $p, q$ with $q \vert p-1$. We also show that certain
twisting deformation of the symmetric group is simple as a Hopf
algebra. On the other hand, we prove that every twisting
deformation of a nilpotent group is semisolvable. We conclude that
the notions of simplicity and (semi)solvability of a semisimple
Hopf algebra are not determined by its tensor category of
representations.
\end{abstract}

\maketitle
\section{Introduction and Main Results}

Let $G$ be a finite group. The character table of $G$ provides
substantial information about the group $G$ itself like, for
instance, normal subgroups and their orders, the center $Z(G)$,
the nilpotence or solvability, etc. In particular, whether $G$ is
simple (respectively, solvable) or not can be established by
inspection of its character table \cite{isaacs}. The character
table provides the same information about $G$ as does the
Grothendieck ring of the tensor category $\Rep G$ \cite{CR}. The
above information is thus \textit{a fortiori} determined by $\Rep
G$.

Finite dimensional Hopf algebras with tensor equivalent categories
of representations are obtained from one another by a twisting
deformation. Properties of $H$ invariant under twisting are of
special interest because they depend only on the tensor category
$\Rep H$.

In this paper we give a series of examples showing that the
notions of simplicity and (semi)solvability of a (semisimple) Hopf
algebra are \emph{not} twist invariants; that is, they are not
categorical notions.

In a first part we show that certain twists of the symmetric group
$\mathbb S_{n}$ on $n$ letters are simple as Hopf algebras, for $n
\geq 5$ (Theorem \ref{condicion}). To prove this we give a
necessary condition for a group-like element in the dual of a
twisting of a finite group to be central, in the case where the
twisting is 'lifted' from an abelian subgroup.

Let $p$, $r$ and $q$ be prime numbers such that $q$ divides $p-1$
and $r-1$. In a second part we show that a family of supersolvable
groups $G$ of order $prq^2$ can be deformed through a twist into
nontrivial simple Hopf algebras (Theorem \ref{A simple dim4pq}).
These twists are also lifted from an abelian subgroup of $G$. The
proof relies on the comparison of the (co)representation theory of
the given twistings \cite{eg-reptriang} with that of an extension
\cite{MW}. We also make use of the classification of semisimple
Hopf algebras in dimension $p$ and $pq$ \cite{zhu, masuoka-pp,
eg-pq, gw}.

It is known that a semisimple Hopf algebra of dimension $p^n$ is
always semisolvable \cite{masuoka, MW}. On the other hand, we
prove that if the group $G$ is nilpotent then any twisting of $kG$
is semisolvable.

Our results imply:

\smallbreak \emph{(a) There exists a simple semisimple Hopf
algebra which is neither twist equivalent to a simple group nor to
the dual of a simple group.}

This answers negatively Question 2.3 in \cite{andrusk}.

\smallbreak \emph{(b) There exists a semisimple Hopf algebra of
dimension $p^2q^2$ which is simple as Hopf algebra.}

Therefore the analogue of Burnside's $p^aq^b$-Theorem for finite
groups does not hold for semisimple Hopf algebras. This concerns
an open question raised by S. Montgomery; see \cite[Question
4.17]{andrusk}.

\smallbreak \emph{(c) There exists a nontrivial semisimple Hopf
algebra which is simple in dimension $36$.}

This example gives a negative answer to \cite[Question, pp.
269]{montgomery}. In dimension $<60$ this is the only possible
such Hopf algebra, by  \cite{Na1}.

\smallbreak \emph{(d) There exists a semisimple Hopf algebra which
is a bosonization but not an extension.}

This answers Question 2.13 of \cite{andrusk}. Indeed, the Hopf
algebras in Theorem \ref{A simple dim4pq} can be built up  as a
Majid-Radford biproduct or bosonization.

\smallbreak We show that there are exactly two twistings of groups
of order $60$ that can be simple as Hopf algebras: the twisting of
$\mathbb A_5$ constructed in \cite{_N}, and the (self-dual)
twisting of $D_3 \times D_5$ discussed in Subsection \ref{60}.
This contributes to the problem in \cite[Question 2.4]{andrusk}.

\smallbreak The paper is organized as follows. In Section
\ref{prels} we recall known facts on normal Hopf subalgebras and
the twisting construction. In Section \ref{symm} we prove some
results on dual central group-like elements in twisting of finite
groups, and present the construction for the symmetric groups.
Section \ref{super} contains the construction for a family of
supersolvable groups; at the end of this section we state our
results in dimensions $60$ and $36$. Finally, in Section
\ref{nilpotentes} we discuss twisting of nilpotent groups.

\smallbreak Along this paper we shall work over an algebraically
closed base field $k$ of characteristic zero. The notation for
Hopf algebras is standard: $\Delta$, $\epsilon$, $S$, denote the
comultiplication, counit and antipode, respectively.

\noindent \textbf{Acknowledgement.} The authors thank N.
Andruskiewitsch for stimulating discussions and encouragement.

\section{Preliminaries}\label{prels}

We discuss in this section basic facts on normal Hopf subalgebras
and twisting deformations.

\subsection{Normal Hopf subalgebras} Let $A$ be a finite dimensional Hopf algebra over $k$.
The (left) adjoint action of $A$ on itself is defined by $\
^hx=\sum h_1x(Sh_2)$,  $x, h\in A$. A Hopf subalgebra $K \subseteq
A$ is called \emph{normal} if it is stable under the adjoint
action; $A$ is called \emph{simple} if it contains no proper
normal Hopf subalgebras.
 Dualizing the notion of normal Hopf subalgebra, we
obtain the notion of conormal quotient Hopf algebra. The notion of
simplicity is self-dual; that is, $A$ is simple if and only if
$A^*$ is.

Let $K \subseteq A$ be a normal Hopf subalgebra. Then $B = A/AK^+$
is a conormal quotient Hopf algebra and the sequence of Hopf
algebra maps $k \longrightarrow K \longrightarrow A
\longrightarrow B \longrightarrow k$ is an exact sequence of Hopf
algebras. In this case we shall say that  $A$ is an
\emph{extension} of  $B$ by  $K$. As a Hopf algebra, $A$ is
isomorphic to a bicrossed product $A \cong K {}^{\tau}\#_{\sigma}
B$.

\smallbreak We recall the following results for future use.

\begin{prop}\cite[Corollary 1.4.3]{Na1}. Let $K \subseteq A$ be a normal Hopf
subalgebra. Suppose that $\dim K$ is the least prime number
dividing $\dim A$. Then $K$ is central in $A$. \qed
\end{prop}

\begin{corol}\label{cental-dimprimo} Let $\pi: A \to B$ be a conormal quotient Hopf algebra.
Suppose that $\dim B$ is the least prime number dividing $\dim A$.
Then $G(B^*)\subset Z(A^*)\cap G(A^*)$. \qed
\end{corol}

\begin{ejem} Let $G$ be a finite group. The normal Hopf subalgebras
of $kG$ are of the form $kH$ where $H$ is a normal subgroup of
$G$. In particular, $kG$ is simple as a Hopf algebra if and only
if $G$ is a simple group. In addition, if $kG$ is simple, then it
possesses no nontrivial quotient Hopf algebra. Therefore $k^G$
contains no proper Hopf subalgebra. \end{ejem}

\subsection{Twisting}
Let $A$ be a finite dimensional Hopf algebra. The category $\Rep
A$ of its finite dimensional representations is a finite
\emph{tensor category} with tensor product given by the diagonal
action of $A$ and unit object $k$.

Finite tensor categories of the form $\Rep A$ are characterized,
using tannakian reconstruction arguments, as those possessing a
fiber functor with values in the category of vector spaces over
$k$. The forgetful functor $\Rep A \to \vect_k$ is a fiber functor
and other fiber functors correspond to \emph{twisting} the
comultiplication of $A$ in the following sense.

\begin{defin} \cite{drinfeld}. A \emph{twist} in $A$ is an invertible $J\in A\otimes A$
satisfying:
\begin{align}
\label{deltaj}(\Delta\otimes id)(J)(J\otimes1)=&(id\otimes\Delta)(J)(1\otimes J),\\
(\varepsilon\otimes id)(J)=&(id\otimes\varepsilon)(J)=1.
\end{align}
\end{defin}
If $J\in A\otimes A$ is a twist, $(A^{J}, m, \Delta^J,
\varepsilon, S^J)$ is a Hopf algebra with $A^J = A$,
$\Delta^J(h)=J^{-1}\Delta(h)J$, and $S^J(h)= v^{-1}S(h)v$, $h\in
A$, $v=m\circ (S\otimes id)(J)$.

The Hopf algebras $A$ y $A'$ are called \emph{twist equivalent} if
$A'\cong A^J$. It is known that $A$ and $A'$ are twist equivalent
if and only if $\Rep A \cong \Rep A'$ as tensor categories
\cite{scha}. Therefore, properties like (quasi)trian\-gularity,
semisimplicity or the structure of the Grothendieck ring are
preserved under twisting deformation.

\begin{obs}\label{morfis tw}
Let $\pi:A \to B$ be a Hopf algebra map and let $J\in A\otimes A$
be a twist. Then $(\pi\otimes\pi)(J)$ is a twist for $B$ and $\pi:
A^J\to B^{(\pi\otimes \pi)(J)}$ is a Hopf algebra map.
\end{obs}

Note that if  $J\in K\otimes K$ is a twist for the Hopf subalgebra
$K\subset A$, then $J\in A\otimes A$ is also a twist for $A$. We
shall say that such $J$ is \emph{lifted} from the Hopf subalgebra
$K$ \cite{ev, vainerman}.

If $A = kG$ is a group algebra, with $G$ a finite abelian group,
the (gauge) equivalence classes of twists for $A$ are in bijective
correspondence with the group $H^2(G, k^*)$ \cite[Proposition
3]{Mv}. The twist $J$ corresponding to the cocycle $\omega: G
\times G \to k^*$ is given by
\begin{equation}\label{twsit abel}
J=\sum_{\alpha,\beta\in
\widehat{G}}\omega(\alpha,\beta)e_{\alpha}\otimes e_{\beta},
\end{equation}where $e_{\chi}=\frac{1}{|G|}\sum_{h\in G}\chi(h^{-1})h$, $\chi \in \widehat G$, is a basis
of orthogonal central idempotents of  $kG$.

\smallbreak Twists in finite groups have been completely
classified \cite{Mv, eg-triangular}. Every twist in $kG$ is lifted
from a \emph{minimal} subgroup $H \subseteq G$ for $J$; that is,
the components of $J_{21}^{-1}J$ span $kH$. Gauge equivalence
classes of twists are classified by classes of pairs $(H,\omega)$,
where $H$ is some solvable subgroup whose order is a square, and
$\omega \in H^2(H, k^*)$ is a non-degenerate 2-cocycle on $H$.

When $H$ is abelian, the twist corresponding to $(H, \omega)$ is
given by \eqref{twsit abel}.

\smallbreak Recall that an element $g \in H$ is called
$\omega$-regular if $\omega(g, h) = \omega(h, g)$, for all $h \in
Z_H(g)$ (this definition depends only on the class of $g$ under
conjugation). Then the cocycle $\omega \in H^2(H, k^*)$ is
non-degenerate if and only if $\{1\}$ is the only $\omega$-regular
class in $H$.

In particular, if $H$ is abelian, the cocycle $\omega$ is
non-degenerate if and only if the skew-symmetric bilinear form
$\omega_{21}^{-1}\omega: H \times H \to k^*$, $(g, h) \mapsto
\omega(g, h)\omega^{-1}(h, g)$, is non-degenerate.

The following lemma follows from \cite{eg-reptriang}. See also
\cite[Lemma 2.11]{_N}.

\begin{lem}\label{no coconmuta}
Let  $J\in kG\otimes kG$ be the twist associated to the pair $(H,
\omega)$, where $H$ is the minimal subgroup of $J$. Then $(kG)^J$
is cocommutative if and only if $H \trianglelefteq G$, $H$ is
abelian and $\omega$ is $\ad G$-invariant in $H^2(H, k^*)$. \qed
\end{lem}

In particular, every nonsymmetric twist lifted from an abelian
subgroup of a simple nonabelian group, or from an abelian subgroup
not containing normal subgroups of $G$, gives rise to a
noncocommutative Hopf algebra.

\begin{ejem}\cite{_N}.\label{ejemplo A_5}
Let $\mathbb A_n$ be the alternating group in  $n$ elements.
Consider the subgroup $H\cong \mathbb{Z}_2\times\mathbb{Z}_2$
generated by $a=(12)(34)$ and $b=(13)(24)$. Let $\omega$ be a
2-cocycle whose cohomology class is nontrivial. Since $\omega$ is
not symmetric, by Lemma \ref{no coconmuta}, $A = (k\mathbb A_n)^J$
is a simple noncommutative and noncocommutative Hopf algebra, for
all $n \geq 5$.

As an algebra, $(k\mathbb A_5)^J \cong k \times M_3(k)^{(2)}
\times M_4(k)\times M_5(k)$. \end{ejem}

Recall that there is a one to one correspondence between quotient
Hopf algebras of $H$ and hereditary subrings of the Grothendieck
ring $K_0(H)$ \cite{NR}. The subring corresponding to the quotient
$H \to \overline H$ is $K_0(\overline H) \subseteq K_0(H)$.

\begin{prop}\label{nik}\cite{_N}.
Let $G$ be a finite simple group and let $J \in kG \otimes kG$ be
a twist. Then $(kG)^J$ is simple as a Hopf algebra. \qed
\end{prop}

Indeed, the result is true under the weaker assumption that $J$ is
a pseudo-twist in $kG$. We shall see that the converse of
Proposition \ref{nik} is not true.

\section{Simple deformations of the symmetric groups}\label{symm}

Let $(A,R)$ be a finite dimensional quasitriangular  Hopf algebra.
The map $f: G(A^*)\to G(A)$, given by $f(\eta)=
R^{(1)}\eta(R^{(2)})$, where $R = R^{(1)}\otimes R^{(2)}$, is an
antihomomorphism of groups. In addition, $\eta \in G(A^*)$ is
central in $A^*$ if and only if $f(\eta)$ is central in $A$
\cite[Proposition 3]{Rad}.

\smallbreak In what follows we shall consider a finite group $G$
and an abelian subgroup $H \subseteq G$. Let $A= (kG)^J$, where
$J\in kG\otimes kG$ is a twist lifted from $H$, written in the
form \eqref{twsit abel}. Observe that $G(A^*) = \widehat G$.

\begin{teor}\label{condicion_centro}
 Suppose that $Z(G) = 1$. Let $\eta \in \widehat G$. Then $\eta\in G(A^*)\cap Z(A^*)$ if
and only if $\eta\vert_H$ is $\omega$-regular.

If $\omega$ is non-degenerate, $\eta\in G(A^*)\cap Z(A^*)$ if and
only if $\eta\vert_H = 1$.
\end{teor}

Here, $\eta\vert_H$ is the restriction of $\eta$ to $H$.

\begin{proof} The Hopf algebra $A$ is triangular with $R$-matrix
$$R = J_{21}^{-1}J = \sum_{\alpha,\beta\in
\widehat{H}}\omega(\alpha,\beta)\omega^{-1}(\beta,\alpha)e_{\alpha}\otimes
e_{\beta}.$$ Since $Z(G) = 1$, we have $Z(A)\cap G(A)=1$. That is,
$\eta\in G(A^*)\cap Z(A^*)$ if and only if $f(\eta)=1$. Let $\eta
\in G(A^*)$. By the orthogonality relations, $\eta(e_{\chi}) =
\delta_{\chi,\eta\vert_H}$, for all $\chi \in \widehat H$.
Therefore,
\begin{align*}
f(\eta) = \sum_{\alpha,\beta\in
\widehat{H}}\omega(\alpha,\beta)\omega^{-1}(\beta,\alpha)e_{\alpha}\eta(e_{\beta})
= \sum_{\alpha\in
\widehat{H}}\omega(\alpha,\eta\vert_H)\omega^{-1}(\eta\vert_H,\alpha)e_{\alpha}.
\end{align*}
Then $f(\eta) = 1$ if and only if $\omega(\chi,\eta\vert_H)
\omega^{-1}(\eta\vert_H,\chi)=1$, $\forall \chi \in \widehat{H}$.
\end{proof}

\begin{corol}\label{ordendeh} Suppose $\omega$ is non-degenerate. Then the order of $H$ divides $[A^*: G(A^*)\cap Z(A^*)]$.
\end{corol}

\begin{proof}In this case, $\eta\vert_H = 1$, for all $\eta \in G(A^*)\cap Z(A^*)$, by Theorem \ref{condicion_centro}.
Then the projection $A \to k(G(A^*)\cap Z(A^*))^*$ restricts
trivially to $H$. The corollary follows from \cite{NZ}.
\end{proof}

Let $\pi: A\to kF$ be a Hopf algebra quotient with $F$ an abelian
group. Then the group $\widehat F \cong F$ can be identified with
a subgroup of $G(A^*)$.

\begin{teor}\label{symm-omega} Suppose $Z(G) = 1$.  Let $\pi: A \to kF$ be a quotient Hopf algebra, where $F \cong \mathbb
Z_p$ and  $p$ is the least prime dividing $|G|$. Then $\pi$ is
normal if and only if $\mu\vert_H$ is $\omega$-regular, for all
$\mu \in \widehat F$. Assume $\omega$ is non-degenerate. Then
$\pi$ is normal if and only if $\mu\vert_H = 1$, for all $\mu \in
\widehat F$. \end{teor}

\begin{proof}By Theorem \ref{condicion_centro}, condition $\omega(\chi,
\mu\vert_H)= \omega(\mu\vert_H, \chi)$, for all $\chi \in
\widehat{H}$, $\mu \in \widehat F$, is equivalent to $\widehat F
\subset Z(A^*)\cap G(A^*)$. In view of Proposition
\ref{cental-dimprimo} this is equivalent to $\pi$ being  normal.
\end{proof}

Let $\pi: \s_n\to \mathbb{Z}_2$ be the only nontrivial
epimorphism. So that $\ker \pi = \mathbb A_n$. Let $H \subseteq
\s_n$ be an abelian subgroup, and let $A=(k\s_n)^J$ be a twisting
with $J\in kH \otimes kH$ given by \eqref{twsit abel}. As noted in
Remark \ref{morfis tw}, $\pi:A\to k\mathbb{Z}_2$ is a Hopf algebra
map.

Consider the sign representation $\sigma : \s_n \to k^*$.

\begin{corol}\label{condicion_prop}$\pi: A \to k\mathbb{Z}_2$
is normal if and only if $\sigma\vert_H$ is $\omega$-regular.

Assume $\omega$ is non-degenerate. Then $\pi$ is normal if and
only if $H \subseteq \mathbb A_n$. \qed
\end{corol}

Let $n \geq 4$. Consider the abelian subgroup $H = <t_1, t_2>
\cong \mathbb{Z}_2 \times \mathbb Z_2$ of $\s_n$, generated by the
transpositions $t_1 = (12)$, $t_2 = (34)$.

We have $\widehat{H} = <a_1, a_2>$, where $a_i(t_j)= 1$ if $i \neq
j$, and $a_i(t_i)= -1$.

Let $\omega$ be the unique nontrivial cocycle on $\widehat H$ (up
to coboundaries). Then $\omega$ is non-degenerate. Let $J \in
kH^{\otimes 2}$ be the corresponding twist.

\begin{teor}\label{condicion}Suppose $n\geq 5$.
Then $(k\s_{n})^J$ is simple.
\end{teor}

\begin{proof} In this case $\pi : \s_n \to \mathbb Z_2$ is the
only nontrivial quotient of $\s_n$. Since twisting preserves
Grothendieck rings, $A = (k\s_n)^J$ has also a unique nontrivial
quotient $\pi : A \to k \mathbb Z_2$. Since $H \nsubseteq \mathbb
A_n$, by Corollary \ref{condicion_prop}, $\pi$ is not normal.
Therefore $(k\s_n)^J$ is simple, as claimed.
\end{proof}

\begin{obs} Let $A = (k\s_n)^J$ as in Theorem \ref{condicion}.
Then $A^*$ is simple and it is not a twisting of any group.
Otherwise $A^*$ and $A$ would be both quasitriangular and simple,
hence $G(A) \cong G(A^*)$ \cite[Proposition 4]{rad-mq}. But this
is not the case, since $H \subseteq G(A)$ and $|G(A^*)| = 2$.
\end{obs}

Another family of examples arises from the construction in
\cite{Bc}. Let $t_i$ denote the transposition $(2i-1,2i) \in
\s_{2n}$, $1 \leq i \leq n$. Consider the abelian subgroup
$H=<t_i, \, 1\leq i\leq n> \cong (\mathbb{Z}_2)^n$ of $\s_{2n}$.

We have $\widehat{H} = <a_i: 1\leq i \leq n>$, where $a_i(t_j)= 1$
if $i \neq j$, and $a_i(t_i)= -1$. Note that $\sigma\vert_H =
a_1a_2\cdots a_n$.

Consider the bicharacter $\omega:\widehat{H}\times\widehat{H}\to
k^*$, $\omega(a_i,a_j)= -1$, $i< j$, $\omega(a_i,a_j)= 1$, $i\geq
j$. This example does not fulfill the condition in Proposition
\ref{condicion_prop} for $n\geq 2$ even and $a = a_1$, since we
have $\omega(\sigma\vert_H,a_1) = (-1)^{n-1}=-1$, while
$\omega(a_1,\sigma\vert_H) = 1$. Let $J$ be  the corresponding
twist.

\begin{teor}\label{condicion-2}Suppose $n\geq 3$, $n$ even.
Then $(k\s_{2n})^J$ is simple.
\end{teor}

\begin{proof} Again in this case $\pi : \s_n \to \mathbb Z_2$ is the
only nontrivial quotient of $\s_n$. So that $A = (k\s_n)^J$ has
also a unique nontrivial quotient $\pi : A \to k \mathbb Z_2$. By
Theorem \ref{condicion} $\pi$ is not normal. Thus $(k\s_n)^J$ is
simple.
\end{proof}

\begin{obs} Let us see that for any twist $J \in k\s_4 \otimes k\s_4$, $(k\s_4)^J$ is not simple;
see \cite[Chapter 6]{Na1}. We know from \cite{Mv, eg-triangular}
that the minimal subgroup $H$ for $J$ is some solvable subgroup
whose order is a square admitting a non-degenerate 2-cocycle. Then
for $k\s_4$ a twist must be lifted from a subgroup $H$ of order
$4$ isomorphic to $\mathbb Z_2 \times \mathbb Z_2$.

\begin{afir} Let $A=(k\s_4)^J$, where $J$ is the twist lifted from a
subgroup $H\cong\mathbb Z_2\times\mathbb Z_2$ which is not normal
and $\omega\neq 1$. Then $G(A)\cong D_4$. \end{afir}

\begin{proof} Let $D \cong D_4$ be the  Sylow 2-subgroup of $\s_4$
containing $H$. Then there are Hopf algebra inclusions $kH\cong
(kH)^J\hookrightarrow (kD)^J\hookrightarrow(k\s_4)^J$.

Since $(kD)^J$ is not commutative and of dimension $8$,
$(kD)^J\cong kD$ as Hopf algebras \cite{TY}. Thus $8$ divides
$|G((k\s_4)^J)|$. As $H$ contains no subgroup which is normal in
$\s_4$, by Lemma \ref{no coconmuta}, $(k\s_4)^J$ is not
cocommutative. Then $|G((k\s_4)^J)|=8$ and $G(A)\cong D\cong D_4$.
\end{proof}

Let the quotients $B = \s_4/K\cong \s_3$,  and $\zeta:(k
\s_4)^J\to (k\s_3)^{\zeta(J)}\cong k\s_3$. We claim that $\zeta$
is normal.

Indeed, we have $\dim A^{co B}=4$. Then $\dim A^{co B}\cap kG(A) =
\dim kG(A)^{co B}$ $=1, 2, 4$, since $|kG(A)|=8$. Also, $\dim
kG(A)^{co B}\dim\zeta(kG(A)) = 8$. If $\dim kG(A)^{co B} = 1,2$,
then $\dim\zeta(kG(A))=8,4$, which is impossible. Then $\dim A^{co
B}\cap kG(A)=4$; that is, $A^{co B}\subset kG(A)$. Hence $A^{co
B}$ is a normal Hopf subalgebra.
\end{obs}

\section{Deformations of a family of supersolvable groups}\label{super}

Let $p$, $q$ and $r$ be prime numbers such that $q$ divides $p-1$
and $r-1$. Let $G_1 = \mathbb Z_p \rtimes \mathbb Z_q$ and $G_2 =
\mathbb Z_r \rtimes \mathbb Z_q$ be the only nonabelian groups of
orders $pq$ and $rq$, respectively. Let $G = G_1 \times G_2$ and
let $\mathbb Z_q \times \mathbb Z_q \cong H \subseteq G$ a
subgroup of order $q^2$. In particular, $G$ is supersolvable and
$Z(G) = 1$.

\medbreak Let $1\neq \omega \in H^2(\widehat H, k^*)$,  $J\in kG
\otimes kG$  the twist lifted from $H$ corresponding to  $\omega$.
Let also $A = (kG)^J$. Note that the cocycle $\omega$ is
nondegenerate. Also, $A$ is a nontrivial Hopf algebra of dimension
$prq^2$.

\begin{lem}\label{alg-coalg}
$A \cong k^{(q^2)} \times \underset{p+r-2 \, \text{copies}}{M_q(k)
\times \dots \times M_q(k)} \times
\underset{\frac{(p-1)(r-1)}{q^2} \, \text{copies}}{M_{q^2}(k)
\times \dots \times M_{q^2}(k)}$ as an algebra.
\end{lem}

\begin{proof} As an algebra, $A = kG \simeq kG_1 \otimes kG_2$. \end{proof}

The coalgebra structure of $A$ follows from the result in
\cite{eg-reptriang} on representations of cotriangular semisimple
Hopf algebras. In particular:

\begin{lem}\label{gp-4} $G(A) \cong H$ is of order $q^2$. \qed \end{lem}

\begin{proof} Let, for all $g \in G$, $H_g = H \cap gHg^{-1}$, and
let $\widetilde\omega$ be the 2-cocycle on $H_g$ given by
$\widetilde\omega(x, y) = \omega^{-1}(g^{-1}xg, g^{-1}yg)
\omega(x, y)$. By \cite{eg-reptriang} the irreducible
representations of $A^*$ are classified by pairs $(\bar g, X)$,
where $\bar g \in H\backslash G/H$ is a double coset modulo $H$
and $X$ is an irreducible representation of the twisted group
algebra $k_{\widetilde\omega}H_g$. The dimension of the
representation $W_{\bar g, X}$ corresponding to $(\bar g, X)$ is
$\dim W_{\bar g, X} = [H: H_g] \dim X$.

Note that $\widetilde \omega$ is trivial on $H_1 = H$. Thus $\dim
W_{\bar 1, X} = 1$, for all possible choices of $X$, giving $|H|$
distinct one-dimensional representations.

Also, $\dim W_{\bar g, X} = 1$ if and only if $H_g = H$ and $\dim
X = 1$. That is, if and only if $g \in N_G(H)$ and $\dim X = 1$.
In our example $N_G(H) = H$. So if $\dim W_{\bar g, X} = 1$, then
$\bar g = \bar 1$ and there are no more one-dimensional
representations. \end{proof}

\begin{lem}\label{abel-ext} Suppose $A$ is not simple and let $K \subseteq A$
be a proper normal Hopf subalgebra. Then $\dim K = pq$ or $rq$.

Moreover, $K$ is necessarily commutative and not cocommutative and
there is an exact sequence of one of the forms $$k \to k^{G_1} \to
A \to kG_2 \to k, \quad k \to k^{G_2} \to A \to kG_1 \to k.$$
\end{lem}

\begin{proof} Let $B = A/AK^+$. There is a Hopf algebra inclusion $B^* \subseteq A^*$ and $B^*$ is
normal in $A^*$. Since $Z(G)=1$, we have  $G(A) \cap Z(A) = 1$. In
particular, $\dim K \neq q$ \cite{zhu}. Suppose $\dim B^* = q$.
Then $B^* = kG(B^*)$ and $G(B^*) \subseteq G(A^*) \cap Z(A^*)$.
Since $\omega$ is non-degenerate, Corollary \ref{ordendeh} implies
that $q^2 \vert [A^*: B^*] = prq$, which is impossible. Hence
$\dim B \neq q$.

If $\dim B = q^2p, q^2r, q^2$, then $\dim K = r, p, pr$ and $K$ is
a group algebra or a dual group algebra \cite{zhu, masuoka-pp,
eg-pq, gw}. Thus $A$ has group-like elements of order $p$ or $r$,
contradicting $|G(A)|=q^2$. Then $\dim K = pq, rq$. Since $|G(A)|
= q^2$, $K\cong k^{G_1}$ or $k^{G_2}$. Similarly, these are the
only possibilities for $B^*$.
\end{proof}

\begin{obs} The results in \cite{eg-reptriang} actually imply that
$A^* \cong A$ as algebras. Once we have established that $A$ is
simple, then we will see in Subsection \ref{selfdual} that
$A^{*\cop} \cong A$ as Hopf algebras. \end{obs}

\begin{teor}\label{A simple dim4pq}
$A$ is simple as a Hopf algebra. \end{teor}

\begin{proof} Suppose not.  We shall compute the dimensions of the irreducible $A$-modules
to get a contradiction. By Lemma \ref{abel-ext}, without loss of
generality, we may assume that $A$ is a bicrossed product $A \cong
k^{G_1}{}^{\tau}\#_{\sigma}kG_2$. By \cite{MW} the irreducible
$A$-modules are classified by pairs $(x, U)$, where $x$ is a
representative of an orbit of the action of $G_2$ in $G_1$ and $U$
is an irreducible  projective representation of the stabilizer
$(G_2)_x$ of $x$. The irreducible module $W_{(x, U)}$
corresponding to $(x, U)$ is the induced module $W_{(x, U)} =
\Ind_{B}^A(kx \otimes U)$, where $B = k^{G_1}\#_{\sigma}(kG_2)_x$.
Therefore $\dim W_{(x, U)} = [G_2: (G_2)_x] \dim U$.

Hence the dimension of an irreducible module cannot be $q^2$. This
contradicts Lemma \ref{alg-coalg}. The theorem is now established.
\end{proof}

\begin{obs} We point out the following consequence of Theorem \ref{A simple dim4pq}.
By Lemmas \ref{alg-coalg} and \ref{gp-4}, $A$ is a biproduct $A =
R \# kH$ in the sense of Majid--Radford, with $\dim R = pr$. This
gives  a nontrivial braided Hopf algebra structure on $R$ over
$H$, in such a way that the corresponding biproduct is
\emph{simple}. This answers \cite[Question 2.13]{andrusk}.
\end{obs}

Letting $p = r$ in Theorem \ref{A simple dim4pq} we obtain:

\begin{teor} Let $p$, $q$, be prime numbers such that $q \vert p-1$. Then there exists a
semisimple Hopf algebra of dimension $p^2q^2$ which is simple as a
Hopf algebra. \qed
\end{teor}

This theorem disproves a conjectured 'quantum version' of
Burnside's $p^aq^b$-Theorem in the context of semisimple Hopf
algebras. It also gives a negative answer to \cite[Question
2.3]{andrusk}.

\subsection{Self-duality}\label{selfdual} If $G$ is finite group
and $J$ is a minimal twist in $G$, then the twisted group algebra
$A = (kG)^J$ satisfies $A\cong A^{*cop}$ \cite{eg-triangular}. We
show next that this is also true under other restrictions on $G$
and $J$.

\begin{prop}
Let $G$ be a finite solvable group and let $J \in kG \otimes kG$
be any twist. Assume $A=(kG)^J$ is simple. Then $A\cong A^{*cop}$.
\end{prop}

By comparing the descriptions of the representation theories of
$G$ and $(kG^J)^*$, we see that this proposition imposes severe
restrictions on the possible groups $G$ satisfying the
assumptions.

\begin{proof} Since $A$ is simple, then $Z(G) = 1$.
Suppose $F \trianglelefteq G$. Then $G$ acts on $k^F$ by the
adjoint action and the smash product $k^F \# kG$ is a Hopf algebra
quotient of $D(G)$ (actually, $D(G)$ corresponds to $F = G$). For
the Hopf algebra $D_F(G) = k^F \# kG$ we have $G(D_F(G))\cap
Z(D_F(G)) = \widehat F \times Z(G) = \widehat F$ and
$$D_F(G)/D_F(G)(k\widehat F)^+ \cong k^{[F, F]} \# kG = D_{[F,
F]}(G).$$ Denote $G = G^{(0)}$, $G^{(i+1)} = [G^{(i)}, G^{(i)}]$,
$i \geq 0$. Since $G$ is solvable, there exists $m \geq 1$ such
that $G^{(m)} = 1$. Iterating the construction above, we get a
sequence
$$D(G) \overset{\pi_1}\to D_{G^{(1)}}(G) \overset{\pi_2}\to D_{G^{(2)}}(G)
\overset{\pi_3}\to \cdots \overset{\pi_m}\to D_{G^{(m)}}(G) =
kG,$$ where every map $\pi_i$ has central kernel. Since $A$ is
twist equivalent to $kG$ then $D(A)$ is twist equivalent to
$D(G)$. In view of Lemma \ref{lema-nil}, we get another sequence
of Hopf algebra maps with central kernels
$$D(A) \overset{\pi_1}\to K_1 \overset{\pi_2}\to K_2
\overset{\pi_3}\to \cdots \overset{\pi_m}\to K_m = (kG)^{J'}.$$

Since the maps $\pi_i$ are normal and $A$ and $A^{*\cop}$ are
simple and nontrivial, we may assume that the composition
$\pi_1\circ \cdots\circ\pi_m:D(A)\to (kG)^{J'}$ is injective when
restricted to $A$ and to $A^{*\cop}$. By dimension, $A\cong
(kG)^{J'} \cong A^{*cop}$.
\end{proof}

\subsection{Twisting deformations of groups of order $60$}\label{60}

As we saw in Example \ref{ejemplo A_5}, there is a nontrivial
simple Hopf algebra of dimension 60, obtained as a twisting of the
alternating group $\mathbb A_5$. Another example arises as a
twisting deformation of the group $D_3\times D_5$, by Theorem
\ref{A simple dim4pq}. For this example, $A \cong A^{*\cop} \cong
k^{(4)} \oplus M_2(k)^{(6)} \oplus M_4(k)^{(2)}$ as algebras; so
$A$ is not a twisting of $k\mathbb A_5$. We shall prove that these
are the only simple Hopf algebras that can arise as twistings of
groups of order $60$.

For the rest of this subsection, $G$ will be a group of order
$60$, $J \in kG \otimes kG$ a twist, and $A = (kG)^J$ a nontrivial
deformation. Since $A$ is not trivial, the minimal subgroup $H$
associated to $J$ must be of order $4$. Then necessarily $H \cong
\mathbb Z_2 \times \mathbb Z_2$ and $J$ corresponds to the only
 nontrivial 2-cocycle on $H$.

We assume first that $G$ is not simple.

\begin{lem}\label{semidirecto} Suppose $A$ is simple. Then $G \cong D_3 \times D_5$.
\end{lem}

\begin{proof} First note that since $G$ is  not simple, then $G$ contains a unique subgroup of
order $5$. Second, the subgroup $H$ cannot be contained in a
normal subgroup $P$ of $G$, since otherwise, the Hopf subalgebra
$kP^J\subset (kG)^J$ would be normal. Similarly, we may assume
that $Z(G) = 1$.

Let $S \trianglelefteq G$ of order $5$. We may assume that  $G' =
G/S$ has a normal subgroup $T$ of order $3$, since otherwise $G$
would contain a normal subgroup of order $20$, thus containing
$H$. Then $N = \pi^{-1}(T) \trianglelefteq G$ is of order $15$; so
$G$ is a semidirect product $N \rtimes H$. Because $Z(G)= 1$,
$G\cong D_3 \times D_5$.
\end{proof}

\begin{teor}\label{sesenta} Let $\vert G\vert = 60$ and let $J \in kG \otimes
kG$ be a twist such that $A = (kG)^J$ is not cocommutative and
simple. Then either \newline \emph{(i)} $G = \mathbb A_5$ and $A$
is isomorphic to the Hopf algebra in Example \ref{ejemplo A_5}; or
\newline \emph{(ii)} $G = D_3 \times D_5$ and $A$ is
isomorphic to the self-dual Hopf algebra in Theorem \ref{A simple
dim4pq}.
\end{teor}

Since the Hopf algebra in (i) is not self-dual, this gives three
simple examples of semisimple Hopf algebras in dimension $60$.

\begin{proof} We use Lemma \ref{semidirecto}. Note that the subgroups of order $4$ in
$G$ are conjugated, and $\vert H^2(\mathbb Z_2 \times \mathbb Z_2,
k^*)\vert = 2$. Hence every pair of nontrivial twists in $G$ gives
rise to isomorphic Hopf algebras.
\end{proof}

\subsection{Twisting deformations of groups of order $36$}\label{36}
The main construction of this section gives an example of a
nontrivial semisimple Hopf algebra $A$ of dimension $36$ that is
simple as a Hopf algebra.  This is thus the smallest semisimple
Hopf algebra which is not semisolvable and the unique simple case
in dimension $36$ \cite{Na1}. The Hopf algebra $A$ is a twisting
of the group $D_3 \times D_3$, and we have $A \cong A^{*\cop}
\cong k^{(4)} \oplus M_2(k)^{(4)} \oplus M_4(k)$ as algebras.

\begin{teor} Let $G$ be a group of order $36$ and let $J \in kG \otimes
kG$ be a twist such that $(kG)^J$ is simple. Then $G \cong D_3
\times D_3$ and $(kG)^J \cong A$.
\end{teor}

\begin{proof} We may assume $Z(G) = 1$. If $J$ is a non-degenerate
twist, then necessarily $G = \mathbb Z_3 \times \mathbb A_4$ or $G
= \mathbb Z_2 \ltimes (\mathbb Z_3 \times \mathbb Z_6)$. This
contradicts the assumption $Z(G) = 1$; hence, we may assume that
$J$ is not minimal. Moreover, the minimal subgroup $H$ of $G$ has
order $4$ or $9$ and it is not cyclic. Also, $H$ is not contained
in any normal subgroup. This leads to consider the case $H \simeq
\mathbb Z_2 \times \mathbb Z_2$  and $G \cong D_3 \times D_3$.
Moreover all subgroups of order $4$ are conjugated and the
conjugation action must preserve nontrivial twists.
\end{proof}

\section{Twisting of nilpotent groups}\label{nilpotentes}
We have shown simple Hopf algebras obtained as twisting of
supersolvable groups. We now prove that this cannot arise from
nilpotent groups.

Let $H$ be a finite dimensional Hopf algebra over $k$.

\begin{defin}  \cite{MW}. A \emph{lower normal series} for $H$ is a series of
proper Hopf subalgebras $H_{n}=k\subset
H_{n-1}\subset\cdots\subset H_1\subset H_0=H$, where $H_{i+1}$ is
normal in $H_i$, for all $i$. The \emph{factors} are the quotients
$\overline{H}_i= H_i/H_iH_{i+1}^+$.

An \emph{upper normal series} is inductively defined as follows.
Let $H_{(0)}=H$. Let $H_i$ be a normal Hopf subalgebra of
$H_{i-1}$ and define $H_{(i)}=H_{(i)}/H_{(i)}H_i^+$. Assume that
$H_n=H_{(n-1)}$, for some positive integer $n$ such that
$H_{(n)}=k$. The \emph{factors} are the Hopf subalgebras $H_i$ of
the quotients $H_{(i)}$.
\end{defin}

Let $G$ be a finite group and let $A=(kG)^J$ be a twisting.

\begin{lem}\label{lema-nil}
Let $Z\subset G$ be a central subgroup. Then $kZ\subset A$ is a
central Hopf subalgebra and $A/A(kZ)^{+}\cong
(kG/Z)^{\overline{J}}$.
\end{lem}

\begin{proof} Since $A=kG$ as algebras,  $kZ$ is central and
$\Delta^J(a)= a\otimes a$. Let $\pi: kG^J\to
k(G/Z)^{\overline{J}}$ be the Hopf algebra map induced by the
 projection $G \to G/Z$. Since $kZ\subset A^{co \pi}$ and
 $\dim A= \dim A^{co\pi} \dim \pi(A)$,  $kZ = A^{co
\pi}$.
\end{proof}

\begin{teor}\label{series-normales} Suppose $G$ is nilpotent. Then
\begin{enumerate}
    \item  $A$ has an upper normal series with factors $k\mathbb Z_p$,
    $p\vert \dim A$,  prime.
    \item $A$ has a lower normal series whose factors are cocommutative.
\end{enumerate} In particular, $A$ is semisolvable in the sense of \cite{MW}.
\end{teor}

\begin{proof} (1). Since $G$ is nilpotent, $Z(G)\neq 1$. Let $Z \subset Z(G)$
be a subgroup of order $p$, $p$ prime. Let $H_1 = kZ \cong
k\mathbb Z_p$ and  $H_{(1)}=A/AH_1^+$. By Lemma \ref{lema-nil},
$H_{(1)}\cong k(G/Z)^{\overline{J}}$. Since $G/Z$ is nilpotent,
 (1) follows by induction on $\vert G\vert$.

(2). Let $H\subset G$ be the subgroup such that $J\in kH\otimes
kH$ is minimal. Since every subgroup of a nilpotent group is
subnormal \cite{Robinson}, we have $H_0=H\vartriangleleft
H_1\vartriangleleft\cdots\vartriangleleft H_n=G$. Then
\begin{equation}\label{serie normal incompleta}
kH_0=kH^J\vartriangleleft
kH_1^J\vartriangleleft\cdots\vartriangleleft kH_n^J=kG^J,
\end{equation}
is part of a lower normal series of $A$ with factors
$k[H_{i+1}/H_i]$, since $J\in kH_i$, for all $i$. Since $H$ is
nilpotent, there is an upper normal series
$H_{(0)}=kH^J\twoheadrightarrow H_{(1)}\cdots\twoheadrightarrow
H_{(s)}=k$.  Moreover,  $kH^J\cong (kH^J)^*$ because  $kH^J$ is
minimal. Thus the dual of this series, that is, $k\hookrightarrow
H_{(1)}^*\hookrightarrow\cdots \hookrightarrow (kH^J)^*\cong
kH^J$, is a lower normal series for $kH^J$ that completes the
series \eqref{serie normal incompleta}.
\end{proof}

\bibliographystyle{amsalpha}

\end{document}